\documentclass[preprint,10pt]{elsarticle}

\makeatletter
\def\ps@pprintTitle{%
 \let\@oddhead\@empty
 \let\@evenhead\@empty
 \def\@oddfoot{}%
 \let\@evenfoot\@oddfoot}
\makeatother

\usepackage{amssymb}
\usepackage[utf8]{inputenc}
\usepackage{amsmath}
\usepackage{amsfonts}
\usepackage{enumitem}  
\usepackage[ruled,vlined,linesnumbered,resetcount]{algorithm2e}
\usepackage{tikz}
\usetikzlibrary{positioning}
\usepackage{caption}
\usepackage{subcaption}

\usepackage{xcolor}
\usepackage{graphicx}
\usepackage[normalem]{ulem}

\usepackage{bm}

\RequirePackage{doi}
\usetikzlibrary{arrows.meta}

\usepackage{minibox}
\usepackage[normalem]{ulem}

\newcommand{\Mourad}[1]{\textcolor{black}{#1}}

\newtheorem{theorem}{Theorem}

\newtheorem{proposition}{Proposition}

\journal{Computer Methods in Applied Mechanics and Engineering}

\newcommand{\myfrac}[2]{\displaystyle{\frac{#1}{#2}}}
\newcommand{\somme}[3]{\displaystyle{\overset{#3}{\underset{#1=#2}{\sum}}}}

\newcommand{\derivee}[2]{\displaystyle{\frac{d #1}{d #2}}}
\newcommand{\deriveepartielle}[2]{\displaystyle{\partial_{#2} #1}}

\newcommand{\vect}[1]{ {\bm{#1}}}

\newcommand{\ctrlvar}{\nu}

\newcommand{\statevecU}{\vect{u}}

\newcommand{\stateU}{u}

\newcommand{\Smode}[2]{\displaystyle{\Phi_{#1}^{#2}}}

\newcommand{\Tmode}[2]{{\alpha_{#1}^{#2}}}

\newcommand{\Sbase}[1]{\Phi^{#1}}

\newcommand{\Xbasis}[1]{{\Phi_{_{#1}}}}

\newcommand{\Tbasis}[1]{{V_{_{#1}}}}

\newcommand{\integOmega}[1]{\displaystyle{\int_{\Omega}#1 \,dx}}

\newcommand{\nbrModes}{q}
\newcommand{\nbrModesU}{q_{u}}
\newcommand{\nbrModesP}{q_{p}}

\newcommand{\nbrParam}{N_{p}}
\newcommand{\nbrSnap}{N_{s}}

\newcommand{\dimVecorsSnap}{N_x}

\newcommand{\trainedSnapMat}[1]{\bm{U}_{_{\trainingparam{_{#1}}}}}

\newcommand{\mcal}[1]{\mathcal{#1}}
\newcommand{\mbb}[1]{\mathbb{#1}}

\newcommand{\txt}[1]{\textnormal{#1}}

\newcommand{\GrassmanifoldSpace}[1]{\mathcal{G}(\nbrModes, N_{#1})}

\newcommand{\TangsubSpaceGrass}[2]{\mcal{T}_{_{[{\Mat{#1}}_{#2}]}}\GrassmanifoldSpace{}}

\newcommand{\Mat}[1]{\bm{#1}}

\newcommand{\Matrix}[3]{\bm{#1^{#2}_{\bm{#3}}}}
\newcommand{\Matint}[3]{{{#1}^{#3}_{#2}}}
\newcommand{\SmodeU}[2]{\Phi^{#1}_{#2}}
\newcommand{\TmodeU}[2]{\alpha^{#1}_{#2}}

\newcommand{\SbasisU}[1]{\bm{\Phi}_{#1}}

\newcommand{\TbasisU}[1]{\bm{\alpha}_{#1}}

\newcommand{\statevarP}{p}

\newcommand{\trainingparam}[1]{\ctrlvar_{#1}}

\newcommand{\InitVel}[1]{\xi_{_{#1}}}

\newcommand{\statefluctvecU}{\statevecU '}
\newcommand{\statefluctvarP}{p'}

\usepackage[a4paper, total={5.5in, 9.in}]{geometry}

\usepackage{floatrow}
\newfloatcommand{capbtabbox}{table}[][\FBwidth]

\SetKwInput{KwOffline}{Offline}
\SetKwInput{KwOnline}{Online}
\usepackage[export]{adjustbox}
\begin{document}

\begin{frontmatter}


% \author{Name\corref{cor1}\fnref{label2}}
% \address{Address\fnref{label3}}

\author{M. Oulghelou\fnref{label1}}
%\ead{}
%\tnotetext[label1]{}
\fntext[label1]{mourad.oulghelou@univ-lr.fr}
%\cortext[cor1]{}
\author{C. Allery\fnref{label2}}
%\tnotetext[label2]{}
\fntext[label2]{cyrille.allery@univ-lr.fr}
%\cortext[cor2]{}

\address{LaSIE, UMR-7356-CNRS, Universit\'e de La Rochelle P\^ole Science et Technologie,
Avenue Michel Cr\'epeau, 17042 La Rochelle Cedex 1, France.}

%\title{Using a hyper reduced variant of the bi-calibrated interpolation on the Grassmann manifold in conjunction with a genetic algorithm to solve inverse flow problems}

%\title{Near real time flow control using genetic algorithm in conjunction with bi-calibrated interpolation on the Grassmann manifold}

\title{A Riemannian Barycentric Interpolation : Derivation of the Parametric Unsteady Navier-Stokes Reduced Order Model}

%\author{Mourad Oulghelou\corref{cor1}, Cyrille Allery\corref{cor2}}
% \address[a]{mourad.oulghelou@univ-lr.fr}
% \address[b]{cyrille.allery@univ-lr.fr}
% \address{}

\begin{abstract}
A new application of subspaces interpolation for the construction of nonlinear Parametric Reduced Order Models (PROMs) is proposed. \Mourad{This approach} is based upon \Mourad{the Riemannian} geometry of the manifold formed by the quotient of the set of full-rank $N$-by-$q$ matrices by the orthogonal group of dimension $q$. 
By using a set of untrained parametrized Proper Orthogonal Decomposition (POD) subspaces of dimension $q$, the subspace for a new untrained parameter is obtained as the generalized \Mourad{Karcher} barycenter which solution is sought after by solving a simple fixed point problem. Contrary to existing PROM approaches, the proposed barycentric PROM is by construction easier to implement and more flexible with respect to change in parameter values.
To assess the potential of the barycentric PROM, numerical experiments are conducted on the parametric flow past a circular cylinder and the flow in a lid driven cavity when the value of Reynolds number varies. It is demonstrated that the proposed barycentric PROM approach achieves competitive results with considerably reduced computational cost.
\end{abstract}


\end{frontmatter}
\newpage
% %%%%%%%%%%%%%%%%%%%%%%%%%%%%%%%%%%%%%%%%%%%%%%%%%%%%%%%%%%%%%%%%%%%%%%%%%%%%%%%%%%%%%%%%%%%%%%%%%%%%%%%%%%%%%%%%%%%%%%%%%%%%%%%%%%%%%%%%%%%%%%%%%%%
\section{Introduction}
Many physical processes governed by parameterized partial differential equations require the ability to predict in real-time, their behavior with respect to variation in parameter values. Unfortunately, this can't be achieved by the standard discretization techniques such as Finite Element, Finite Volumes or Finite Differences. Therefore, techniques of model reduction must be used. 
The basic objective in model reduction is to turn a large scale high fidelity dynamical system into a low-order initial value \Mourad{problem, able} to accurately reproduce the dynamics of the original system. To achieve such low-order systems, a projection subspace where the high fidelity solutions can be accurately reproduced needs to be determined. The temporal dynamics is thereafter described by the reduced order model (ROM) \Mourad{usually} obtained via Galerkin projection of the original system onto the projection subspace. The method that is most used to determine projection subspaces is the POD (Proper Orthogonal Decomposition) method \cite{Sirovich}. Basically, the POD consists in solving an eigenvalue problem formed \Mourad{from the dynamical} system snapshots picked at different time instants of the physical problem. Only the most energetic first modes corresponding to highest eigenvalues are considered in the construction of the POD basis. The POD approach has been applied in many fields such as micro electro-mechanics  \cite{LIANG2002}, aeroelasticity \cite{Thomas2003, Lucia2004}, structural dynamics \cite{Park2007,ALDMOUR2002,Buljak2011}, damage detection \cite{Galvanetto2008, Shane2011}, modal analysis \cite{Han2003}, chemical reaction dynamics \cite{ Singer2009}, fluid mechanics \cite{Allery-2005, Beghein-Allery-2014}, etc.
A major issue of a ROM based on POD\Mourad{,} arises often in parametric problems such as design and optimal control. In general, a POD ROM is valid only in a small neighborhood of the parameters for which the POD basis was built \cite{AKKARI2014522, Bergmann}. Therefore, suitable ROMs covering the variation in parameters need to be derived. 
There is a growing interest in building fast and accurate parametric reduced order models (PROMs) and many attempts have been made on this subject. Hereafter, generic approaches are reviewed and the key idea of the contribution of the present paper is presented.
\vspace*{0.2cm}
\\
\textit{Approach by a priori reduction} : A possible way to build PROMs is by means of the Proper Generlizaed Decomposition (PGD). This method is an a priori model reduction technique that enables to generate a tensorial form of the solution including space, time and parameters of interest. Originally, the PGD was used to achieve a space-time separation of the solution of nonlinear structural mechanics problems \cite{LADEVEZE20033061, LADEVEZE20101287} and was extended to achieve space-parameter separation of the solution of polymeric systems \cite{Ammar2006153,Mokdad2007}. It was since used in different other fields such as stochastic problems \cite{Nouy2007}, heat transfer \cite{Pruliere2010}, quantum chemistry \cite{Chinesta2008}, fluid dynamics \cite{Dumon2011a,  LEBLOND2014}, etc. In the methodological point of view, the construction of the PGD separated form is progressively performed throughout some well defined enrichment functions. These functions are determined simultaneously by solving for many iterations, the subproblems resulting from the double Galerkin projections \cite{Chinesta2011b}. Even though the time-space or parameter-space PGD demonstrated its efficiency for a wide range of applications, it remains a costly method for real-time applications. %This is due to the need to perform many iterations of the high fidelity problem whenever an update of the separated form is needed.  
%In the context of parametric model reduction, the PGD is usually used in the offline stage to create a vade-mecum \cite{Chinesta2013} of the solution, that can thereafter be used online to predict in real-time the solution over a range of parameter values. 
%However to the best of our knowledge, no such attempt have been made successfully in fluid dynamics problems governed by time dependent Navier-Stokes equations.
\vspace*{0.2cm}
\\
\textit{Approach by using a global subspace} : A straight forward approach previously used in the context of reduced optimal control \cite{Ravindran1997, Tallet-Allery-2016}, consists in extending the snapshots ensemble in some sort to construct a POD basis that hopefully covers the dynamics features in different trained parameters. By this global approach \footnote{The PGD can also be used here in the offline stage to create a vade-mecum \cite{Chinesta2013} of the solution. This vade-mecum can thereafter be used online to predict in real-time the solution over a range of parameter values.}, the ROM coefficients obtained via Galerkin projections are calculated once for all during the offline stage, and in the online stage, only the ROM equations are solved whenever the parameter changes. From the computational point of view, the global approach may suffer from requiring a large number of snapshots and thus a large number of modes for the constructed POD basis. Moreover, the information overload in the POD basis affects the ROM predictions of features that occur in a restricted regime. To overcome this issue, parametric  model reduction via  interpolation are privileged. 
\vspace*{0.2cm}
\\
\textit{Approach by subspaces interpolation} : A basic approach to construct PROMs is by using the practical Riemannian geometry of the Grassamnn manifold \cite{Absil, Edelman1998}. Typically, a widely used approach introcuced first in the context of aeroelasticity is the ITSGM (Interpolation on the Tangent Space of the Grassmann Manifold) method  \cite{Amsallem}. 
Constructing a ITSGM PROM consists of five steps : (1) from the trained parametrized POD bases, chose a POD basis to be the reference point of tangency to Grassmann manifold (2) use the Riemannian Logarithm and map all the subspaces spanned by the trained POD bases to this space; (3) perform linear interpolation in the tangent space; (4) map back the result by the Riamannian Exponential to the Grassmann manifold; (5) perform Gelrerkin projections to obtain the PROM which is a system of ordinary differential equations (ODE). The ITSGM has been successfully applied in ROM adaptation \cite{Amsallem, OulghelouBiCITSGM2019Arxiv} adjoint based optimal control \cite{OULGHELOUAMC2018}, data driven optimal control \cite{oulghelou2020}, etc. In the same spirit of the ITSGM, a selection of methods using the Riemannian geometry of the Grassmann manifold were recently proposed. These methods are presented as : a generalization to the Grassmann manifold of the Inverse Distance Weighted (IDW) method \cite{Rolando_IDW_article}, a generalization to the Grassmann manifold of Neville Aitken's method performed recursively by calculating the geodesic barycenter of two points \cite{rolando_Neville}, a kriging technique where interpolation on the tangent space to the Grassmann manifold are performed with sophisticated weighting factors obtained as solution of an optimization problem \cite{TheseRolando}. By all the aforementioned subspace interpolation methods, the update of the parametric reduced order model with respect to change in parameter values lead to non-negligeable recalculations of Galerkin projections. Hereafter, a new flexible PROM procedure that enables a real-time update of  Galerkin projections with respect to parameter change is introduced.
%
%
%Moreover, it is more flexible with the change of parameter values and easier to implement
%
%In all the data sets considered here, we have never faced issues related to ill-definitions of these tools.
%
%Unfortunately, to our knowledge, there has been no work addressing this issue
\vspace*{0.2cm}
\\
\textit{Proposed approach} : Basically, constructing the ROM coefficients via Galerkin projections and solving the resulting ODE system constitutes the online stage of parametric projection model reduction. The main focus of the present paper is to tackles parametric variations in nonlinear PROMs without performing any further Galerkin projections in the online stage. 
For this sake, a PROM construction approach based upon the practical Riemannian geometry of the quotient manifold\footnote{The quotient manifold $\mbb{R}^{N\times \nbrModes}_*/\mcal{O}(q)$ is formed by the quotient of the set of full-rank $N$-by-$q$ matrices by the orthogonal group of dimension $q$.} $\mbb{R}^{N\times \nbrModes}_*/\mcal{O}(q)$ \cite{Massar2020} is proposed.
Initially, this geometry was applied to interpolate low-rank solutions of the Luyapunov equations resulting from parametric linear input-output reduced order systems \cite{Massart2019}. In the present article, the use of this geometry is extended to the general interpolation framework of parametrized nonlinear reduced order models, where the trained points of interpolation are considered to be the parametrized POD subspaces and the untrained subspace is their generalized Karcher barycenter \cite{Karcher}. This leads to an optimization problem for which the solution is found by solving a simple fixed point problem.  In the remaining of the paper, the new PROM approach is referred to as barycentric PROM. The accuracy of the method is assessed on two flow problems when the Reynolds number value changes. These are the flow past a circular cylinder and the flow in a lid driven cavity. In order to not clutter up the study, the barycentric approach is compared only to the standard ITSGM method. Steps of ITSGM are outlined in \ref{Appendix_ITSGM}.
%
%
%
%
%
%
%
%
%%%%
\vspace*{0.2cm}
\\
\indent The paper is organized as follows: a brief overview of the quotient geometry used to derive the barycentric interpolation approach is given in section 2. Illustration of this approach on the construction of the PROM for generic Navier-Stokes equations is given in section 3. In section 4 the potential of the barycentric interpolation is studied and compared to the existing ITSGM on two parametric problems : the flow past a circular cylinder and the flow in a lid driven cavity where the parameter to vary is the Reynolds number. Conclusions of the paper are drawn in the final section.
%
%
%
%
%
%
%
%
%%%%

%\section{Problem statement}
%	\input{./Tex_Files/problem_setting.tex}

\section{Riemannian barycentric interpolation strategy}
In a differential manifold $\mcal{M}$, the straight line between two points $Y_1$ and $Y_2$ that remain in the manifold is called a geodesic path \cite{Wald}. This path is parametrized by a twice differentiable function $\gamma(t)$ (with $0\leq t\leq 1$) solution of a second order initial value problem. Thus, a geodesic path that starts from $Y_1$ and ends in $Y_2$ is defined by two initial conditions. These are the initial value $\gamma(0)=Y_1$ lying in the manifold $\mcal{M}$ and the initial velocity $\dot{\gamma}(0) = \xi$ lying in $\mcal{T}_{Y_1} \mcal{M}$ the tangent space to the manifold $\mcal{M}$ at the point $Y_1$. 
The link between a manifold and its tangent space at a point is established trough the exponential and logarithmic mappings. The exponential map sends a point from the tangent space to the manifold and the logarithm is its inverse. In general, the exponential of $\xi$ is given by
\begin{equation}\label{Eq:Exp_map}
\txt{Exp}_{Y_1}(\xi) = \gamma(1)
\end{equation}
The quotient geometry in which rely the main contribution of the present paper is rigousorly studied in \cite{Massar2020}. This geometry results in cheap practical formulas for the exponential and logarithmic maps that will be recalled in the next subsection. 
	\subsection{Exponential and logarithmic maps on $\mbb{R}^{N\times \nbrModes}_{*}/\mathcal{O}(\nbrModes)$}
Let $\nbrModes < N$ be two positive integers, the set of all the matrices $\SmodeU{}{}\in\mbb{R}^{N\times \nbrModes}_{*}$ that yield the same subspace  is the quotient manifold defined by
$$ \mbb{R}^{N\times \nbrModes}_{*}/\mathcal{O}(\nbrModes) := \mbb{R}^{N\times \nbrModes}_{*}/ \sim \ \ := \{\SmodeU{}{} \mathcal{O}(\nbrModes) \ \ : \ \ \SmodeU{}{}\in \mbb{R}^{N\times \nbrModes}_{*} \}$$
where $\sim$ denotes the equivalent class in $\mbb{R}^{N\times \nbrModes}_{*}$ such that $\SmodeU{1}{} \sim \SmodeU{2}{}$ if and only if $\SmodeU{2}{} = \SmodeU{1}{} Q $ for some orthogonal matrix $Q\in \mathcal{O}(\nbrModes)$.
Define the following map
\begin{equation}
\begin{matrix}
\pi:&\mbb{R}^{N\times \nbrModes}_{*}  &\longrightarrow  &\mbb{R}^{N\times \nbrModes}_{*}/\mathcal{O}(\nbrModes)  
\vspace*{0.15cm}
\\
&\SmodeU{}{} && \SmodeU{}{} \mathcal{O}(\nbrModes)
\end{matrix}
\end{equation}
The set of horizontal vectors at a point 
%$\SmodeU{}{} \mathcal{O}(\nbrModes)\in\mbb{R}^{N\times \nbrModes}_{*}/\mathcal{O}(\nbrModes)$ 
$\SmodeU{}{}\in\mbb{R}^{N\times \nbrModes}_{*}$
is given by
$$\mcal{H}_{\SmodeU{}{}} =  \{  \Phi (\Phi^T \Phi)^{-1} H + \Phi_{\bot} K \ \ \txt{such that} \ \ H = H^T \in \mbb{R}^{\nbrModes \times \nbrModes}, \Phi ^T \Phi_{\bot} = 0_{_{\nbrModes\times (N-\nbrModes)}}  , K\in \mbb{R}^{(N-\nbrModes) \times \nbrModes} \}$$
Note that for an arbitrary $\xi_{\Phi} \in \mcal{H}_{\Phi}$, the curve $t\mapsto \Phi + t \xi_{\Phi} $ does not necessarily remain full rank. In order to be able to define the exponential map on the quotient manifold, the following set of allowed horizontal vectors is introduced
$$ \mcal{D}_{\Phi} = \{ \xi_{\Phi} \in\mcal{H}_{\Phi} \ \ \txt{such that} \ \ \txt{rank}(\Phi+t\xi_{\Phi}) = \nbrModes \ \ \forall t \in [0,1]\} $$
The expression of the exponential map 
%on $ \mbb{R}^{N\times \nbrModes}_{*}/\mathcal{O}(\nbrModes)$ 
is given by the following theorem \cite{Massar2020}
\begin{theorem}
For all $\xi_{\SmodeU{}{}} \in \mcal{D}_{\SmodeU{}{}}$, the exponential map on $ \mbb{R}^{N\times \nbrModes}_{*}/\mathcal{O}(\nbrModes)$ is given by
\begin{equation}\label{Eq.Exponential_map}
\txt{Exp}_{\pi(\SmodeU{}{})}\xi_{\SmodeU{}{}} = \pi(\SmodeU{}{} + \xi_{\SmodeU{}{}})
\end{equation}
\end{theorem}
The reciprocal map of the exponential is known as the logarithmic map. It is determined  by solving for the horizontal lift $\xi$ the equation $\txt{EXP}_{\pi(\Phi)} \xi = \pi(\Psi)$, where $\Phi,\Psi \in \mbb{R}^{n\times \nbrModes}_{*}$. The practical expression of the logarithmic map is given by the following result \cite{Massar2020}
\begin{proposition}\label{prop.Log}
Let $\Phi,\Psi \in \mbb{R}^{n\times \nbrModes}_{*}$ such that $\Phi^T \Psi$ is nonsingular. Then $\txt{Log}_{\pi(\Phi)} \pi(\Psi)$ is uniquely defined as
\begin{equation}
\txt{Log}_{\pi(\Phi)} \pi(\Psi) = \Psi Q - \Phi, \quad Q = VU^T,
\end{equation}
where $U$ and $V$ are respectively the left and right singular matrices of $\Phi^T \Psi$. i.e, $\Phi^T \Psi \overset{SVD}{=} U \Sigma V^T$. Moreover, the distance between $\pi(\Phi)$ and $\pi(\Psi)$ is given by
$$d(\pi(\Phi), \pi(\Psi)) = ||\Psi Q - \Phi||_F$$
\end{proposition}
	\subsection{Proposed barycentric interpolation strategy}
Consider a set of parametrized POD subspaces $\pi(\Sbase{1}), \pi(\Sbase{2}), \dots, \pi(\Sbase{\nbrParam})$ associated respectively to the trained parameter values $\nu_1, \nu_2, \dots, \nu_{\nbrParam}$. The aim of the following is to predict the subspace $\pi(\widetilde{\SmodeU{}{}})$ associated to a new untrained value $\tilde{\nu}\neq \nu_k$. The idea of the barycentric interpolation approach consists in seeking $\pi(\widetilde{\SmodeU{}{}})$ as the Karcher barycenter \cite{Karcher} of the points $\pi(\Sbase{1}), \pi(\Sbase{2}), \dots, \pi(\Sbase{\nbrParam})$. Mathematically, $\pi(\widetilde{\SmodeU{}{}})$ is defined  as the minimizer of the Riemannian objective function
\begin{equation}\label{Riemannian_minimization_problem}
\mcal{J}(\pi(\SmodeU{}{})) = \myfrac{1}{2} \,
\somme{k}{1}{\nbrParam}\omega_k(\tilde{\nu})\, d^2(\pi(\SmodeU{k}{}) , \pi(\SmodeU{}{}))
\end{equation}
where $\{\omega_k : \tilde{\nu} \to \omega_k(\tilde{\nu}), \ \ k=1\dots,\nbrParam \}$ is a set of interpolation functions verifying 
\begin{align*}
\somme{k}{1}{\nbrParam}\omega_k(\tilde{\nu}) = 1 \quad\txt{and}\quad
\omega_k(\nu_h) = \delta_{kh}
\end{align*}
A possible choice of interpolation functions can be Lagrangian functions, inverse distance weighted functions, radial basis functions, etc. Finding $\SmodeU{}{}$ requires to solve the optimization problem \eqref{Riemannian_minimization_problem}. To this end, a local minimum can be found as the root of the gradient of the functional $\mcal{J}$. By definition, the Riemannian gradient $\nabla^{\pi(\SmodeU{}{})} $ of the squared distance $d^2(\pi(\Phi) , \pi(\Psi))$ with respect to a reference point $\pi(\SmodeU{}{})$ is equal to $-2 \, \txt{Log}_{\pi(\SmodeU{}{})} \pi(\Psi)$. This yields to 
$$ \nabla^{\pi(\SmodeU{}{})} \mcal{J}(\pi(\widetilde{\SmodeU{}{}})) = - \somme{k}{1}{\nbrParam}\omega_k(\tilde{\nu}) \, \txt{Log}_{\pi(\SmodeU{}{})} \pi(\SmodeU{k}{}) $$ 
By substituting the expression of the Riemnannian logarithm given in Proposition \ref{prop.Log}, The gradient of $\mcal{J}$ writes 
$$ \nabla^{\pi(\SmodeU{}{})} \mcal{J}(\pi(\widetilde{\SmodeU{}{}})) = \somme{k}{1}{\nbrParam}\omega_k(\tilde{\nu}) \left( \widetilde{\SmodeU{}{}} - \SmodeU{k}{} \widetilde{Q}_k \right) $$
where $\widetilde{Q}_k = V_k U_k^T$ is obtained from the decomposition $\widetilde{\Phi}^T \SmodeU{k}{} \overset{\txt{SVD}}{=} U_k \Sigma_k V_k^T$.
Given that the sum of the weights $\omega_k(\tilde{\nu})$ is equal to 1, when $\nabla^{\SmodeU{}{}} \mcal{J}(\pi(\widetilde{\SmodeU{}{}}))$ vanishes, the following fixed point problem yields
$$\widetilde{\SmodeU{}{}} = \somme{k}{1}{\nbrParam}\omega_k(\tilde{\nu}) \SmodeU{k}{} \widetilde{Q}_k
$$
%
% is assumed to be nonsingular.
The implementation of the fixed point strategy is outlined in algorithm \ref{Alg:barycenter_SPsD}.
\vspace*{0.2cm}
\\
\begin{algorithm}[H]
Give a value of $\tilde{\nu}$ (chosen by the user) and calculate the weights $\omega_1(\tilde{\nu}), \dots, \omega_{\nbrParam}(\tilde{\nu})$\\
Set $\widetilde{\SmodeU{}{}}$ arbitrary, for example choose a point $\SmodeU{k}{}$ from the sampling\\
\While{$ || \nabla^{\pi(\SmodeU{}{})} \mcal{J}(\pi(\widetilde{\SmodeU{}{}})) ||_F >\varepsilon$}{
\For{$k\in\{1,\dots,\nbrParam\}$}{
Evaluate the matrix $\widetilde{Q}^k = V_k U_k^T$ where $\widetilde{\Phi}^T \SmodeU{k}{} \overset{\tiny{\txt{SVD}}}{=} U_k \Sigma_k V_k^T$
}
Calculate the barycenter : $\widetilde{\SmodeU{}{}} = \somme{k}{1}{\nbrParam} \omega_k(\tilde{\nu}) \SmodeU{k}{} \widetilde{Q}^k$\\
%Evaluate the functional gradient $\nabla^{\SmodeU{}{}} \mcal{J}(\SmodeU{}{})  =  - \somme{k}{1}{\nbrParam}\omega_k \txt{Log}_{\pi(\SmodeU{}{})} \pi(\SmodeU{k}{})$\\
%Assign the current barycenter as a reference point : $\SmodeU{ref}{} \leftarrow \SmodeU{}{}$
}
\caption{Barycentric interpolation of subspaces in $\mbb{R}^{n\times \nbrModes}_*/\mcal{O}(\nbrModes)$}
\label{Alg:barycenter_SPsD}
\end{algorithm}

\section{Illustration of Barycentric PROM on generic Navier-Stokes equations}
Consider the Newtonian flow governed by the Navier-Stokes equations
\begin{equation}\label{EQ.Navier_Stokes}
\begin{cases}
\deriveepartielle{\statevecU}{t} - \nu \Delta \statevecU + \statevecU\cdot\nabla \statevecU + \frac{1}{\rho}\nabla \statevarP  = \bm{f}  & \txt{in } \Omega\times]0,T]\\
\nabla\cdot\statevecU = 0 & \txt{in } \Omega\times]0,T]\\
\statevecU = \bm{g} & \txt{on } \partial\Omega \times]0,T]\\
\statevecU(0) = \statevecU_0 & \txt{in } \Omega
\end{cases}
\end{equation}
where $\statevecU$ is the velocity of initial value $\statevecU_0$, $p$ the pressure, $\bm{f}$ the external body forces , $\nu$ the kinematic viscosity and $\rho$ the fluid density. The computational domain is an open bounded connected subset of $\mbb{R}^d$ with $d=2,3$. On the boundary $\partial\Omega$, the velocity is set to the time independent function $\bm{g}$.
%, and the force per surface unit exerted at each point on $\Gamma_{N}$ is equal to $f_N$.
By varying the kinematic viscosity value, we are interested in predicting the solutions $\statevecU$ of Navier-Stokes equations by using the proposed barycentric interpolation. 
\vspace*{0.2cm}
\\
For $\nu_1, \dots, \nu_{\nbrParam}$ a set of values of $\nu$, assume that the discrete solutions of the high fidelity equations \eqref{EQ.Navier_Stokes} in the partition $t_0<t_1<\dots<t_{\nbrSnap} $ of the computational time interval $[0,T]$ is available. That is, for each $\nu_k$, we have the following snapshots matrix 
\begin{equation*}
\trainedSnapMat{k} = \begin{bmatrix}
\statevecU(t_{_1}, x_1, \trainingparam{k}) & \statevecU(t_{_2}, x_1, \trainingparam{k}) &  \\
&&\\
\vdots & \ddots &  \\
&&\\
\statevecU(t_{_1}, x_{\dimVecorsSnap}, \trainingparam{k}) &  & \statevecU(t_{_{\nbrSnap}}, x_{\dimVecorsSnap}, \trainingparam{k})
\end{bmatrix}
\end{equation*}
where $\dimVecorsSnap$ is the number of degrees of freedom. It is assumed that $\dimVecorsSnap$ exceeds the number of snapshots $\nbrSnap$ by several orders of magnitude, $\dimVecorsSnap\gg \nbrSnap$ and that $\trainedSnapMat{k}$ are rank-$\nbrModes$ matrices, $\nbrModes\leq\nbrSnap$. 
In order to obtain the reduced order model, the solution $\statevecU$ is decomposed as
\begin{equation*}
\statevecU(t,x, \nu_k) = \overline{\statevecU}(x) + \statefluctvecU(t,x_i, \nu_k) 
\end{equation*} 
where $\overline{\statevecU}$ represents the mean field and $\statefluctvecU$ the fluctuations. For each value $\nu_k$, a POD basis of dimension $\nbrModes$ is constructed on the fluctuating part such that
\begin{equation*}\label{EQ : subspace approx POD u and p}
\statefluctvecU(t,x, \nu_k) \approx  \somme{j}{1}{\nbrModes}  \SmodeU{k}{j}(x)  \TmodeU{k}{j}(t)
\end{equation*}
For a new value $\tilde{\ctrlvar} \neq \ctrlvar_k$, the solution by the proposed barycentric interpolation approach writes as follows
%By using the proposed Riemannian barycentric interpolation, the solution for $\tilde{\ctrlvar} \neq \ctrlvar_k$ is expressed by the relation
%
\begin{equation}\label{EQ.general_approx_u_and_p}
\statevecU(t,x, \tilde{\ctrlvar})  \approx \overline{\statevecU}(x) + \somme{k}{1}{\nbrParam} \somme{j,s}{1}{\nbrModes}\omega_k(\tilde{\ctrlvar})   \widetilde{Q}^{k}_{sj}  \SmodeU{k}{s}(x) \widetilde{\alpha}_j (t) 
\end{equation}
Note that in the approximations above, the weights $\omega_k(\tilde{\ctrlvar})$ are predetermined and the orthogonal matrices $\widetilde{Q}_k$ are the outcome of algorithm \ref{Alg:barycenter_SPsD}. It remains now to \Mourad{determine the vector $\widetilde{\TbasisU{}}(t) = \{\widetilde{\alpha}_1 (t), \widetilde{\alpha}_2 (t),\dots, \widetilde{\alpha}_{\nbrModes} (t) \}$} describing the temporal dynamics by solving the PROM which is a system of ordinary differential equations. The barycentric PROM is constructed by performing Galerkin projections of the momentum equation in \eqref{EQ.general_approx_u_and_p} onto the interpolated basis. This yields to
\begin{equation}\label{EQ.ROM_Navier_Stokes}
\Matrix{M}{}{} \derivee{\widetilde{\TbasisU{}}}{t} + \nu \Matrix{R}{}{} \widetilde{\TbasisU{}} + \Matrix{\overline{C}}{}{} \widetilde{\TbasisU{}} +  \somme{e}{1}{\nbrModes}\widetilde{\alpha}_e \Matrix{C}{}{e}  \widetilde{\TbasisU{}}  = \bm{\widetilde{F}}
\end{equation}
with 
\begin{equation*}\label{Eq.ROM_update_U}
\begin{aligned}
&\Matrix{M}{}{} = \somme{k,h}{1}{\nbrParam}  \omega_k(\tilde{\nu}) \omega_h(\tilde{\nu})  Q^{k^T} \Matint{M}{}{hk} Q^{h},
\quad\txt{such that}\quad 
\Matint{M}{ij}{hk} =  \integOmega{ \SmodeU{k}{j} \SmodeU{h}{i} },
\\
&\Matrix{R}{}{}  = \somme{k,h}{1}{\nbrParam}  \omega_k(\tilde{\nu}) \omega_h(\tilde{\nu})  Q^{k^T} \Matint{R}{}{hk} Q^{h},
\quad\txt{such that}\quad 
\Matint{R}{ij}{hk} =  \integOmega{\nabla \SmodeU{k}{j} : \nabla\SmodeU{h}{i}},
\\
&\Matrix{\overline{C}}{}{}  = \somme{k,h}{1}{\nbrParam} \omega_k(\tilde{\nu}) \omega_h(\tilde{\nu})  Q^{k^T} \Matint{\overline{C}}{}{hk} Q^{h},
\quad\txt{such that}\quad 
\Matint{\overline{C}}{ij}{hk} =   \integOmega{ (\overline{\statevecU}\cdot\nabla) \SmodeU{k}{j}\SmodeU{h}{i}  } + \integOmega{ (\SmodeU{k}{j}\cdot\nabla) \overline{\statevecU} \SmodeU{h}{i}},
\\
&\Matrix{C}{}{e} = \somme{k,h,n}{1}{\nbrParam}   \omega_{n}(\tilde{\nu}) \omega_k(\tilde{\nu}) \omega_h(\tilde{\nu}) \somme{s}{1}{\nbrModes} Q^{n}_{se}   Q^{k^T} \Matint{C}{s}{hkn} Q^{h}%, \quad e = 1,\dots,\nbrModes
\quad\txt{such that}\quad 
\Matint{C}{sij}{hkn} =   \integOmega{(\SmodeU{k}{j}\cdot\nabla) \SmodeU{n}{s}\SmodeU{h}{k}  },
\\
&\bm{\widetilde{F}} = \somme{k}{1}{\nbrParam}  \omega_k(\tilde{\nu})  Q^{k^T} F^{k} ,
\quad\txt{such that}\quad 
F^k_i=   \integOmega{ \bm{f} \, \SmodeU{k}{i}} - \tilde{\nu} \integOmega{\nabla \overline{\statevecU} : \nabla \SmodeU{k}{i}} - \integOmega{ (\overline{\statevecU}\cdot\nabla) \overline{\statevecU}\SmodeU{k}{i}  }
\end{aligned}
\end{equation*}
Note that in the above formulas, $\Matrix{M}{}{}$, $\Matrix{R}{}{}$, $\Matrix{\overline{C}}{}{}$ and $\Matrix{C}{}{e}$ are $q$-by-$q$ matrices and $\bm{\widetilde{F}}$ is a vector of length $q$. The superscripts $h$, $k$ and $n$ refer to the indices of the untrained values of the parameter $\nu$ and the subscripts $i$, $j$ and $s$ to the indices of the POD modes. Obviously, the most expensive step in the construction of a ROM is the calculation of its coefficients obtained via Galerkin projections. This is particularly a bottleneck in parametric model reduction where these costly calculations are performed online whenever the parameter $\ctrlvar$ changes. By its natural construction, the proposed barycentric PROM \eqref{EQ.ROM_Navier_Stokes} provides a solution to this issue. The principal advantage of the barycentric PROM lies in the possibility to transfer all the intense calculations resulting from Galerkin projections to the offline stage. This leads to efficient update throughout  cheap manipulations of $\nbrModes \times \nbrModes$ matrices. For instance, the computational complexity to update the reduced matrix $\Matrix{M}{}{} $ by the barycentric PROM is $2\nbrParam^2\mcal{O}(\nbrModes^3) + (\nbrParam^2-1)\mcal{O}(\nbrModes^2)$, while it is $\mcal{O}(\dimVecorsSnap^2 \nbrModes)+\mcal{O}(\dimVecorsSnap \nbrModes^2)$ by the ITSGM PROM, where $\nbrModes$ and $N_p$ are smaller than $ \dimVecorsSnap$ by several orders of magnitude. 
%In practice, the potential of the barycentric PROM in terms computational time saving appears greatly when the number of degrees of freedom $N_x$ is high. 
The steps of the proposed barycentric Navier-Stokes PROM are summarized in algorithm \ref{Alg:PROM_NS}.
\vspace*{0.2cm}
\\
\begin{algorithm}[H]
\begin{itemize}
\item[]
\begin{itemize}
\item[\textit{\textbf{Offline :}}]
\item[\textit{step 1 :}] Let $\nu_1, \dots, \nu_{\nbrParam}$ be a set of parameters. For each $\nu_k$, calculate and store the POD basis $\SbasisU{k}$.
\item[\textit{step 2 :}] For $h,k, n\in\{1,\dots,\nbrParam\}$ and $s\in\{1,\dots,\nbrModes\}$ assemble and store the reduced matrices $\Matint{M}{}{hk}$, $\Matint{R}{}{hk}$, $\Matint{\overline{C}}{}{hk}$, $F^{k}$ and $\Matint{C}{s}{hkn}$ 
%\item[\textit{step 4 :}] For $h,k,n\in\{1,\dots,\nbrParam\}$ assemble and store the tensors $\Matint{C}{[hkn]}{s}$,  $s = 1,\dots,\nbrModes $.
%
%
%\item[\textit{step 3 :}] For each $k_0\in\{1,\dots,\nbrParam\}$, calculate and store the orthogonal matrices $Q_{{[k]}}$, $k = 1,\dots,\nbrParam$, given by
%$$Q_{{[k]}} = Z_{u_{[k]}} W^T_{u_{[k]}} \quad\txt{with}\quad \bm{\singvalU{k_0}{}}\SbasisU{k_0}^T \SbasisU{k}\bm{\singvalU{k}{}} = W_{u_{[k]}} \Theta_{u_{[k]}} Z^T_{u_{[k]}}$$
%
%
\item[\textit{\textbf{Online :}}]
\item[\textit{step 1 :}] Give a value of $\tilde{\nu}$ (chosen by the user) 
\item[\textit{step 2 :}] Chose the weights $\omega_1(\tilde{\nu}),\dots,\omega_{\nbrParam}(\tilde{\nu})$. For instance, Lagrange, RBF, IDW, etc%such that $\omega_{k_0} = 1 - \underset{k\neq k_0}{\somme{k}{1}{\nbrParam}} \omega_k$
\item[\textit{step 3 :}] Perform algorithm  \ref{Alg:barycenter_SPsD} and obtain the orthogonal matrices $\widetilde{Q}^1, \widetilde{Q}^2,\dots, \widetilde{Q}^{\nbrParam} $
\item[\textit{step 4 :}] Update the matrices $\Matrix{M}{}{}$, $\Matrix{R}{}{}$, $\Matrix{\overline{C}}{}{}$, $\bm{F}$ and $\Matrix{C}{}{e}$ for $e \in \{1,\dots,\nbrModes\}$
%
%
%\item[\textit{step 5 :}] Update pressure ROM coefficients $\Matrix{M}{p}{}$, $\Matrix{R}{p}{}$, $\Matrix{\overline{C}}{p}{}$, $\Matrix{K}{p}{}$ and $\bm{F^p}$
%
%
%\item[\textit{step 5 :}] Update the tensor $\Matrix{C}{e}{}$ for $e \in \{1,\dots,\nbrModes\}$
%
%
\item[\textit{step 5 :}] Solve the parametric reduced order model \eqref{EQ.ROM_Navier_Stokes} %and obtain $\TbasisU{}{}$ 
\end{itemize}
\end{itemize}
\caption{Navier-Stokes PROM}
\label{Alg:PROM_NS}
\end{algorithm}
\section{Numerical tests}
\subsection{Discretization of the high fidelity problem}
To state the discretization approach used in practical implementation of the Navier-Stokes equations, we restrict to the case where the density $\rho=1$ and the external forces $\bm{f}=0$. Consider a time step $\delta t>0$, the approximate solutions $\bm{u}^n$ and  $p^n$ are respectively computed by solving at each step $n$, the following problem 
\begin{equation}
\begin{cases}
\statevecU^n - \nu \delta t \Delta \statevecU^n + \myfrac{3 \delta t}{2} \statevecU^{n-1}\cdot\nabla \statevecU^{n-1} - \myfrac{ \delta t}{2} \statevecU^{n-2}\cdot\nabla \statevecU^{n-2} + \delta t\nabla \statevarP^{n}  = \statevecU^{n-1}
\\
\nabla\cdot \statevecU^n = 0
\end{cases}
\end{equation}   
Note that in this scheme, the first order backward Euler method is used to discretize the temporal derivative, implicit scheme for the diffusion term and Adams-Bashfort \cite{Yinnia2010} scheme to deal with the nonlinear convective term. On the other hand,  the Finite Element method is used for the spatial discretization. In order to guarantee the inf-sup stabilization condition, the Taylor-Hood Finite Element pair $\mbb{P}_2/\mbb{P}_1$ is chosen.
The full discretization of the problem in hand gives rise to a sequence of linear symmetric saddle-point systems of the form
\begin{equation}\label{saddle-point_problem}
\begin{cases}
A \statevecU^n -  \delta t B^T \statevarP^{n}  = \delta t  \mcal{F}(\statevecU^{n-1}, \statevecU^{n-2})
\\
B \statevecU^n = 0
\end{cases}
\end{equation} 
%This obviously leads to cheaper Stokes-like problems that can be efficiently implemented in the parallel framework. 
A common problem of saddle-point systems is that they can be poorly conditioned. Thus, care must be taken in their numerical implementation. That is why the strategy based on the iterative augmented Lagrangian method is adopted. More details about this approach can be found in \cite{book-Glowinski}. 
In the following, the parameter of interpolation is the Reynolds number given by 
$$Re = U D/\nu$$ 
where $U$ is the characteristic velocity, $D$ the characteristic length and $\nu$ the kinematic viscosity. The variation of the Reynolds number is carried out by fixing $D$ and $U$ and changing the value of the kinematic viscosity $\nu$.
\vspace*{0.2cm}
\\
Consider a bunch of trained values $\nu_1, \nu_2,\dots,\nu_{\nbrParam} $. The saddle point systems \eqref{saddle-point_problem} are numerically solved by using Fenics \cite{FenicsBook} and the POD bases are constructed on the fluctuating parts of the  velocity and pressure solutions. To this end, the variables $\statevecU$ and $p$ are decomposed as follows
\begin{equation*}
\begin{cases}
\statevecU(t,x, \nu) = \overline{\statevecU}(x) + \statefluctvecU(t,x, \nu),  \hspace{0.5cm}\txt{ in } \Omega \\
p(t,x, \nu) = \overline{p}(x) + \statefluctvarP(t,x, \nu),  \hspace{0.5cm}\txt{ in } \Omega \\
%\overline{\statevecU} = U,  \hspace{0.5cm}\txt{ on } \Gamma_{inflow}\\
\end{cases}
\end{equation*}
The mean parts $\overline{\statevecU}$ and $\overline{p}$ are chosen such as
\begin{equation*}
\overline{\vect{\stateU}} = \myfrac{1}{\nbrParam \nbrSnap} \somme{k}{1}{\nbrParam}\somme{j}{1}{\nbrSnap} \statevecU(t_j, x, \nu_k)
\hspace*{2cm}
\overline{p} = \myfrac{1}{\nbrParam \nbrSnap} \somme{k}{1}{\nbrParam}\somme{j}{1}{\nbrSnap} \statevarP(t_j, x, \nu_k)
\end{equation*}
where $\nbrParam$ is the number of trained parameter values and $\nbrSnap$ the number of snapshots. In order the assess the accuracy of the PROM, the $L^2$ percentage of errors at a time $t$ and the corresponding mean relative error are used. If $\bm{u}_{\txt{ref}}$ and $\bm{u}_{\txt{approx}}$ are respectively the reference high fidelity solution and its approximation obtained by a PROM, the error at a time $t$ and its mean counterpart over a time frame $[t_0, t_1]$ write respectively  as follows
$$ \mcal{E}_u(t)  = 100\times \left. \left( \integOmega{ |\bm{u}_{\txt{ref}}(t)  - \bm{u}_{\txt{approx}}(t)|^2 }\right)^{\frac{1}{2}} \middle/ \left(\integOmega{ |\bm{u}_{\txt{ref}}(t)|^2 }\right)^{\frac{1}{2}} \right.$$

$$ \bar{\mcal{E}}_u  = 100\times \left. \left( \int_{t_1}^{t_2} \integOmega{ |\bm{u}_{\txt{ref}}(t)  - \bm{u}_{\txt{approx}}(t)|^2 }\, dt\right)^{\frac{1}{2}} \middle/ \left(\int_{t_1}^{t_2} \integOmega{ |\bm{u}_{\txt{ref}}(t)|^2 }\, dt\right)^{\frac{1}{2}} \right.$$
%The contribution ratio of the $k^{\txt{th}}$ mode is given by
%$$\txt{RIC}^{k} = \somme{i}{1}{k}\lambda_i / \somme{i}{1}{\nbrSnap}\lambda_i$$
%where $\lambda_i$ are the POD eigenvalues
%
%
%
%
%\footnote{The variation of the Reynolds number is carried out through the variation of the kinematic viscosity $\nu$. The inlet velocity $U$ in this case is kept constant.}
%
%
%
%
The aim in the following examples is on one hand, to test the ability of the Barycentric PROM in predicting well and fast the solutions for the new untrained values of interpolation parameter, and on the other hand, to provide a comparative study in terms of accuracy and computational time with respect to the existing approach ITSGM.
We mention that for both Barycentric and ITSGM methods, the subspace interpolation is carried out by selecting the POD bases associated to the nearest three trained Reynolds number values. 
\subsection{Flow past a circular cylinder}
\begin{figure}[hbtp!]
\centering 
\includegraphics[width=0.85\linewidth]{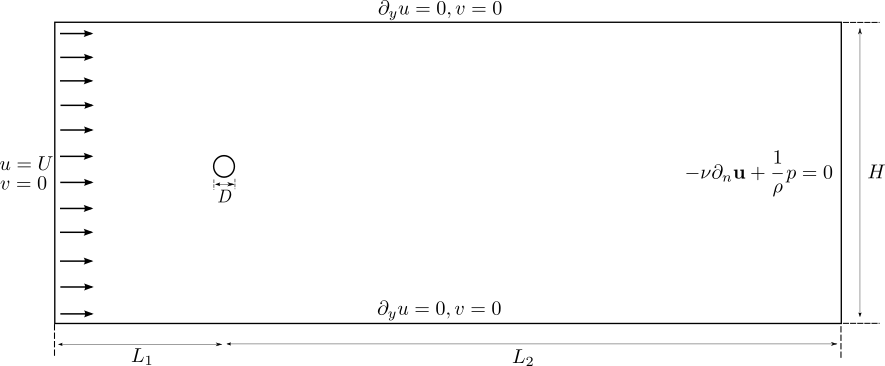}
\caption{Two-dimensional domain and boundary conditions for the problem of flow past a circular cylinder.}
\label{Fig:flow past cylinder}
\end{figure}

\begin{figure}[hbtp!]
\includegraphics[width=0.32\linewidth]{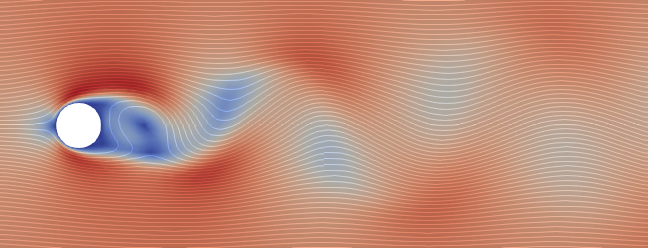}%
\hspace*{0.25cm}\includegraphics[width=0.32\linewidth]{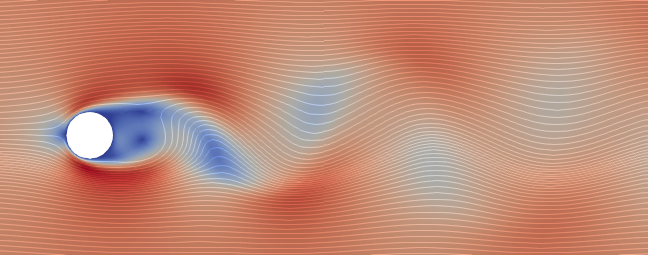}%
\hspace*{0.25cm}\includegraphics[width=0.32\linewidth]{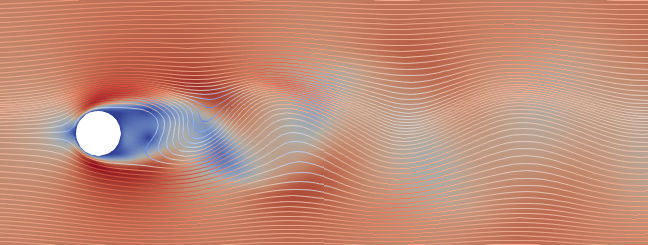}

\caption{Von Kármán vortex streets in three different instants of the established regime of the flow past a cylinder when $Re=160$.}
\label{fig:dynamics_cylinder_flow}
\end{figure}

\begin{figure}[hbtp!]
\centering 
\includegraphics[width=\linewidth]{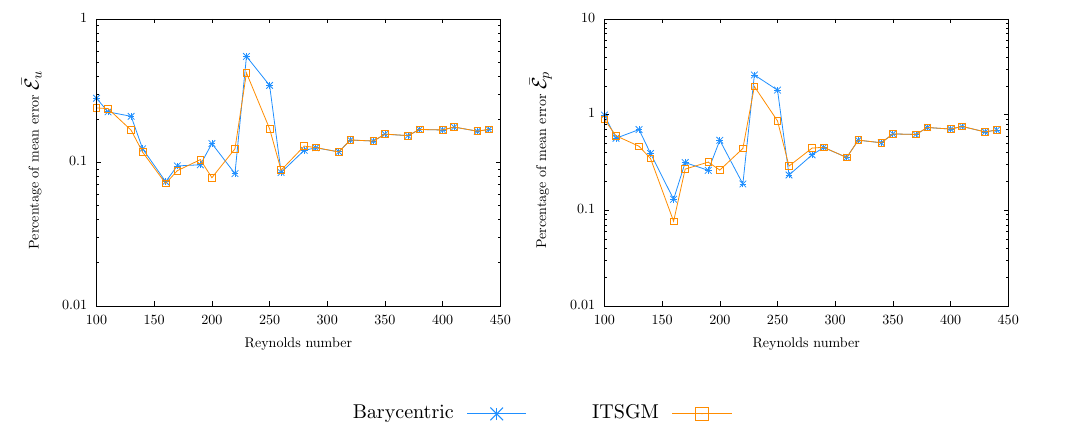}
\caption{Percentage of the mean errors at different untrained Reynolds number values for the proposed barycentric and ITSGM PROMs.}
\label{fig:Errors_Re_range_cylinder}
\end{figure}

\begin{figure}[hbtp!]
\centering 
\includegraphics[width=\linewidth]{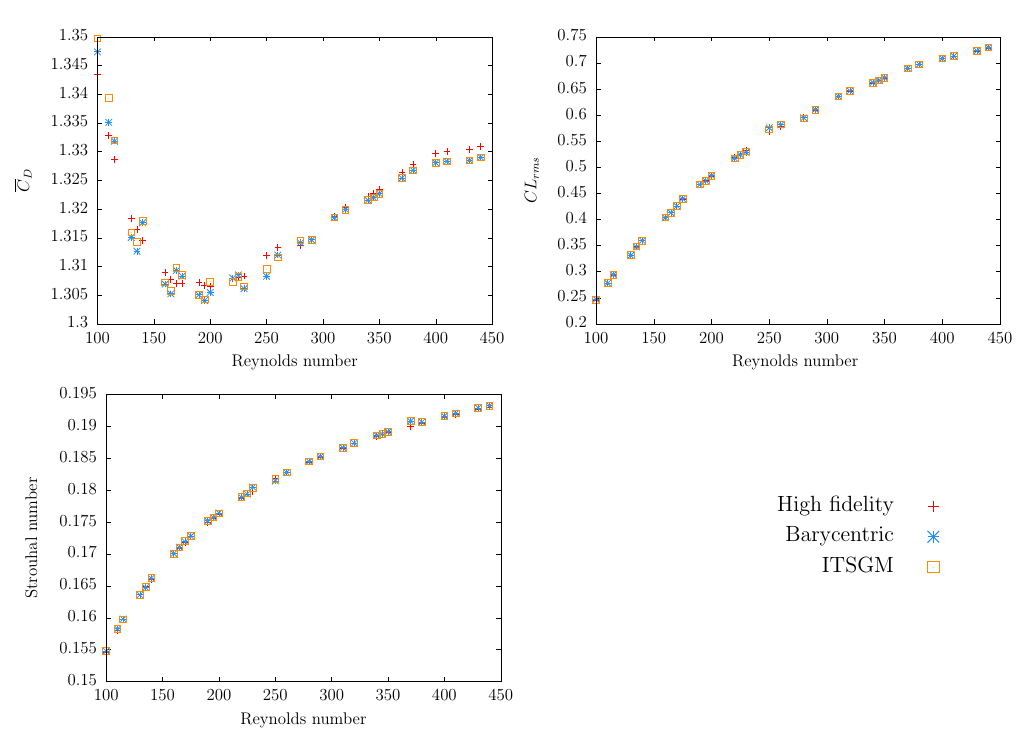}
\caption{Hydrodynamic coefficients obtained with the high fidelity model, the proposed barycentric PROM and the ITSGM PROM.}
\label{fig:hydrodynamic_coeffs}
\end{figure}

Consider the flow in a rectangular domain past a circular cylinder of diameter $D$. The length of the rectangle is $H = 30D$ and its width is $45 D$. The center of the cylinder is located at $L_1 = 10D$ from the left boundary and $H/2$ from the lower boundary. The fluid dynamics of the flow is driven by an inlet velocity $U$ from the left boundary $\Gamma_{\txt{in}}$, and allowed to flow past through the right boundary $\Gamma_{\txt{out}}$. Free slip boundary conditions are imposed on the horizontal edges $\Gamma_{\txt{horiz}}$ and no slip boundary condition on the cylinder’s wall $\Gamma_{\txt{circ}}$. A sketch up of the domain and boundary conditions is given in figure \ref{Fig:flow past cylinder}. In the discrete problem, the time step is set to $\delta t = 0.01$ and a mesh containing $21174$ elements is used for the triangulation of the spatial domain. A representation of the flow dynamics at $Re=160$ in different time instants is given in figure \ref{fig:dynamics_cylinder_flow}. 
It can be seen that a flow pattern forms around the cylinder creating by that the well known periodic lateral von Kármán vortex street.
%
%
% BEGIN FIGURE -----------------------------------------------------------------------------------------------------------------------------
% END FIGURE -------------------------------------------------------------------------------------------------------------------------------
\\
A database of flow solutions is constructed at the trained Reynolds number values $Re^{\txt{tr}}$ ranging from $90$ to $450$, where the jump between two consecutive values is equal to $30$. 
At each value $Re^{\txt{tr}}_k$, $500$ snapshots of the flow covering about $8$ periods of the established regime are selected in order to build the POD bases.
The order of truncation of velocity and pressure variables is respectively set to $\nbrModes_u=10$ and $\nbrModes_p=8$. The aim of this application is to predict the mean drag coefficient, the root mean square lift coefficient and the Strouhal number for a wide range of untrained values of the Reynolds number. The expressions of these hydrodynamics quantities are given in \ref{Appendix_aerodynamics_coeffs}.
To achieve that, the determination by the PROM of the pressure is also required. The PROM \eqref{EQ.ROM_Navier_Stokes} has then to be equipped with an additional equation which is the Galerkin projection of the high fidelity system onto the gradient of pressure POD bases. Accordingly and without loss of generalities, few more steps need to be added to algorithm \ref{Alg:PROM_NS}. 
%In the offline stage includes two more steps : the determination of the pressure POD bases and the calculation of their associated Galerkin projections. The change in the online stage consists in : the barycentric interpolation of pressure bases and the associated reduced matrices update. A description of this PROM is given in appendix \ref{}.
%
%
These are one one hand, the determination of the pressure POD bases and the calculation of their associated Galerkin projection in the offline stage, and on the other hand, the barycentric interpolation with the associated reduced matrices update in the online stage. Details about the used PROM are given in \ref{Appendix_ROM2}.
%As a result, the  the offline stage of algorithm \ref{}  Galerkin projections of the high fidelity system onto the gradient of pressure POD bases is performed. 
In the following, the interpolation of the flow in several untrained Reynolds number values by using the Barycentric and ITSGM approaches is assessed. The untrained values are chosen in the range between $90$ and $450$ with a jump of $10$.
\vspace*{0.2cm}
\\
The accuracy of the barycentric interpolation is first assessed by calculating the percentage of mean $L^2$ errors in the untrained points. Figure \ref{fig:Errors_Re_range_cylinder} reveals that besides the good accuracy of the barycentric PROM, the obtained mean errors have nearly the same behavior with those obtained by the ITSGM PROM. In the worst case scenario, the error of the barycentric PROM is about $0.6\%$ for velocity and $2.6\%$ for the pressure.  In figure \ref{fig:hydrodynamic_coeffs}, the hydrodynamics coefficients obtained by barycentric and ITSGM PROMs are compared to each other and to their high fidelity counterpart. The reader can observe that the root mean square lift coefficient and strouhal number obtained by the barycentric and ITSGM PROMs are nearly identical to their high fidelity counterpart. Slight differences are however observed in the prediction of the mean drag coefficient by both approaches. 
In terms of the computational time of PROMs update with respect to parameter variation, the barycentric PROM is instantaneously updated in $8\times 10^{-3}$ seconds while $7.1$ seconds are spent to update the ITSGM PROM. That is $3$ orders of magnitudes of time saving achieved by the barycentric PROM. %This result supports the statement made about the effectiveness and flexibility of this last approach in real-time applications.

%More specifically
%\begin{figure}[hbtp!]
%\hspace*{-0.6 in}
%\begin{subfigure}{.5\textwidth}
%  \centering
%  \resizebox{\textwidth}{!}{\input{./Figures/Cylinder/mean_Relative_Errors/mean_Relative_Error_velocity}}
%  %\caption{Pourcentage d'erreur de la vitesse}
%  \label{fig:sub1}
%\end{subfigure}%
%\begin{subfigure}{.5\textwidth}
%  \centering
%  \resizebox{\textwidth}{!}{\input{./Figures/Cylinder/mean_Relative_Errors/mean_Relative_Error_pressure}}
%  %\caption{Pourcentage d'erreur de la pression}
%  \label{fig:sub2}
%\end{subfigure}
%\\ 
%\centering
%\hspace*{-1cm}\adjustbox{width= 0.6\linewidth,trim= 1cm 4cm 1cm 4cm}{\input{./Figures/Cylinder/mean_Relative_Errors/LEGEND_errors}}
%\caption{Percentage of the mean errors at different untrained Reynolds number values for the proposed barycentric and ITSGM PROMs.}
%\label{fig:Errors_Re_range_cylinder}
%\end{figure}

%\begin{figure}[hbtp!]
%\hspace*{-2cm}
%\begin{subfigure}{0.5\textwidth}
%\resizebox{\textwidth}{!}{\large \input{./Figures/Cylinder/mean_coeffs/CD_mean}}%
%%\caption{Point $x_1 = (2/16, 13/16)$.}
%\end{subfigure}%
%\begin{subfigure}{0.5\textwidth}
%\resizebox{\textwidth}{!}{\large \input{./Figures/Cylinder/mean_coeffs/Crms}}%
%%\caption{Point $x_2 = (2/16, 2/16)$.}
%\end{subfigure}
%\hspace*{-0.9cm}
%\begin{subfigure}{0.5\textwidth}
%\resizebox{\textwidth}{!}{\large \input{./Figures/Cylinder/mean_coeffs/strouhal}}%
%%\caption{Point $x_3 = (19/20, 19/20)$.}
%\end{subfigure}%
%%\centering
%\hspace*{1.6cm}
%\hspace*{-1cm}\adjustbox{width= 0.55\linewidth,trim= 1cm 4cm 1cm 4cm}{\large \input{./Figures/Cylinder/mean_coeffs/LEGEND_coeffs}}
%\caption{Hydrodynamic coefficients obtained with the high fidelity model, the proposed barycentric PROM and the ITSGM PROM.}
%\label{fig:hydrodynamic_coeffs}
%\end{figure}

\subsection{Lid Driven Cavity flow}
\begin{figure}[hbtp!]
\includegraphics[width=0.32\linewidth]{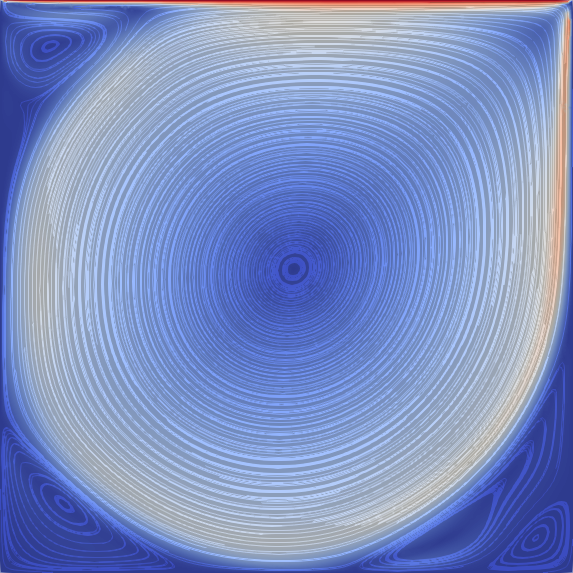}%
\hspace*{0.25cm}\includegraphics[width=0.32\linewidth]{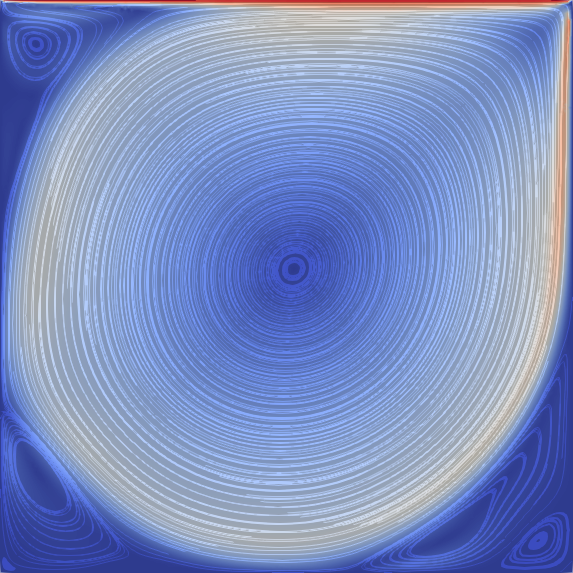}%
\hspace*{0.25cm}\includegraphics[width=0.32\linewidth]{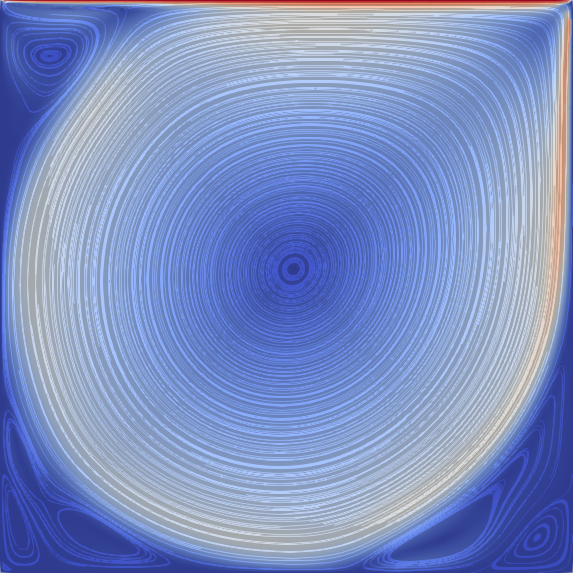}

\caption{Streamlines in three different instants of the established regime of the flow in a lid driven cavity when $Re=9000$.}
\label{fig:dynamics_cavity_flow}
\end{figure}

\begin{figure}[hbtp!]
\centering 
\includegraphics[width=0.55\linewidth]{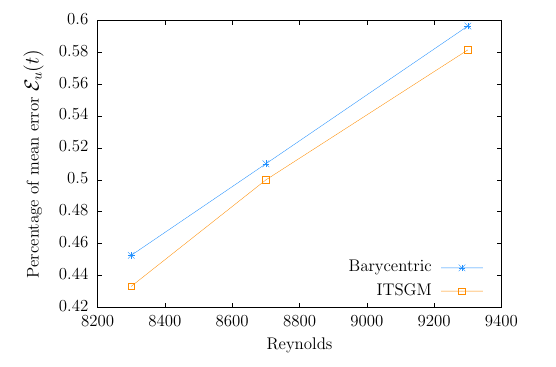}
\caption{Percentage of the mean errors of the solutions by the barycentric and ITSGM PROMs.}
\label{fig:Mean_Errors}
\end{figure}

\begin{figure}[hbtp!]
\centering 
\includegraphics[width=\linewidth]{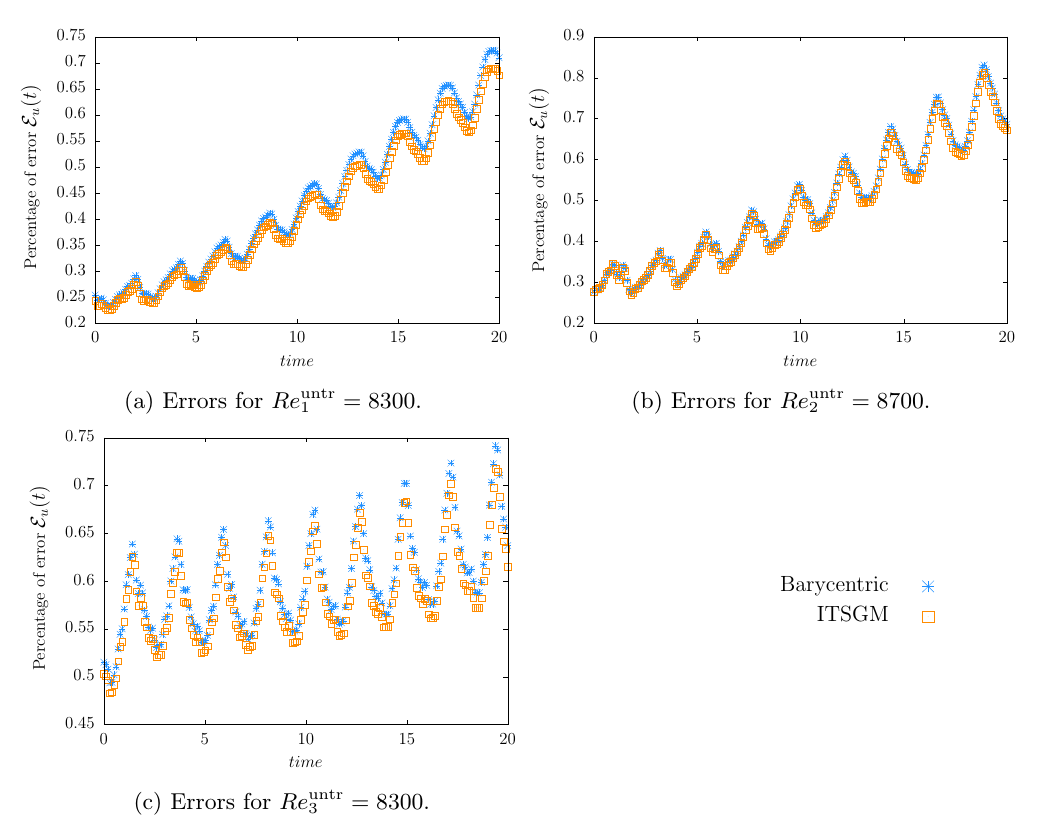}
\caption{Time history of the percentage of error of the solutions obtained by the barycentric and ITSGM PROMs.}
\label{fig:Errors_time}
\end{figure}

\begin{figure}[hbtp!]
\centering 
\includegraphics[width=\linewidth]{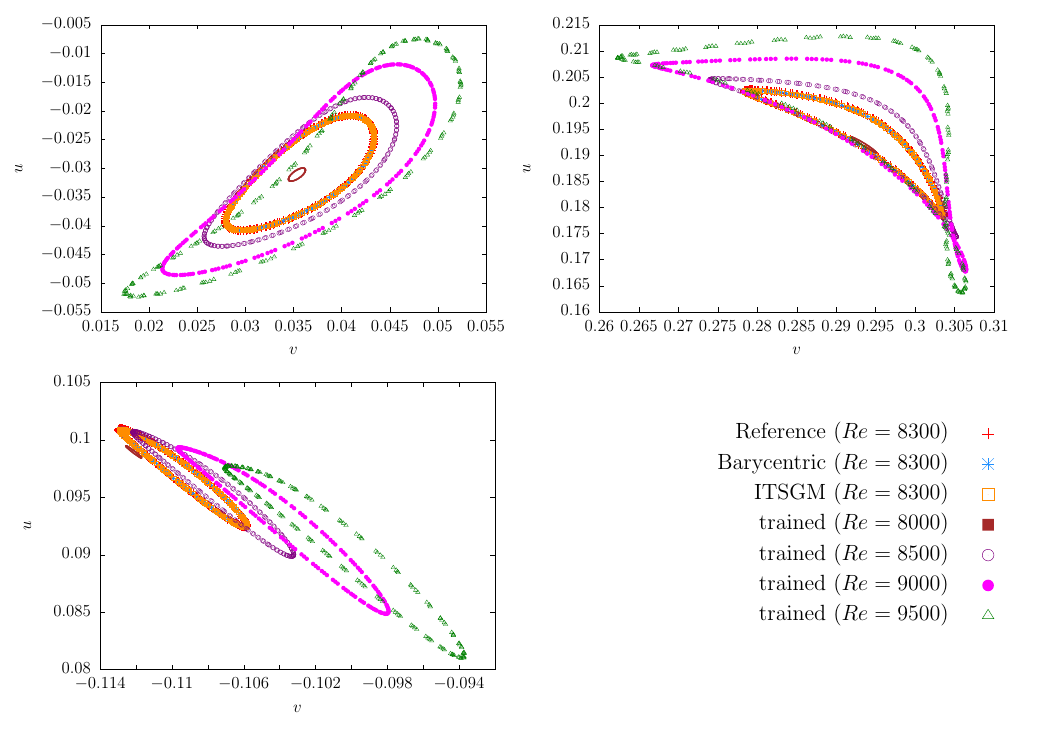}
\caption{Velocity phase portrait at the \Mourad{points $x_1$ (top left), $x_2$ (top right) and $x_3$ (bottom left)} for the untrained value $Re_1^{\txt{untr}} = 8300$. 
}
\label{fig:velocity_at_points_RE8300}
\end{figure}

\begin{figure}[hbtp!]
\centering 
\includegraphics[width=\linewidth]{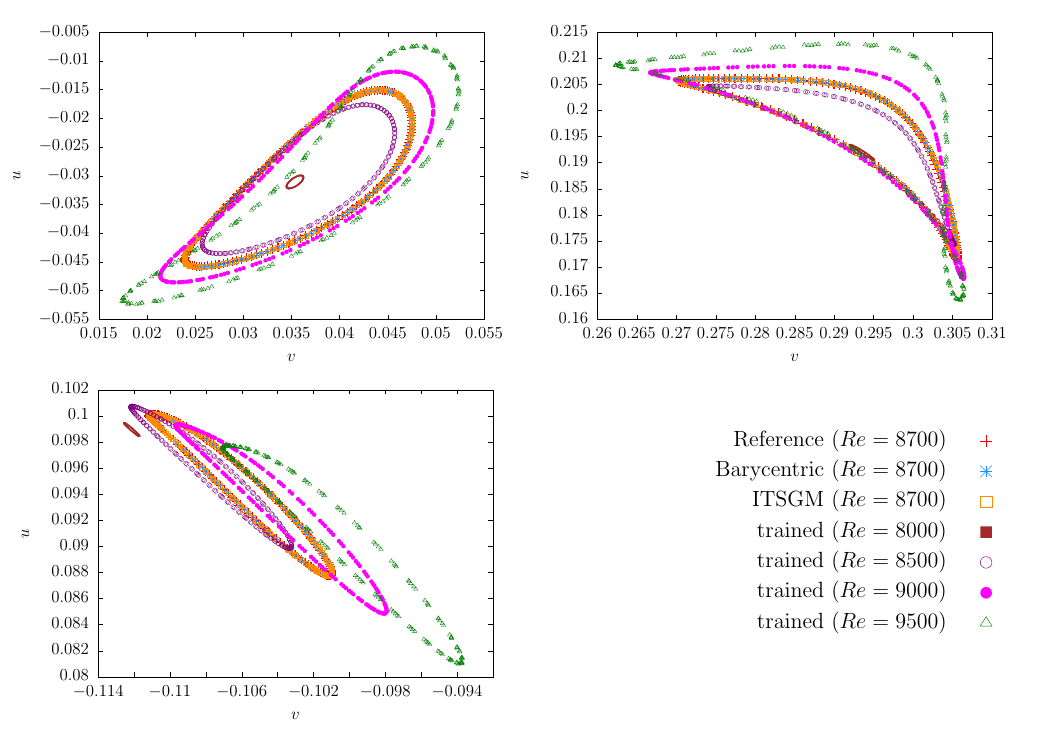}
\caption{Velocity phase portrait at the \Mourad{points $x_1$ (top left), $x_2$ (top right) and $x_3$ (bottom left)} for the untrained value $Re_1^{\txt{untr}} = 8700$. 
}
\label{fig:velocity_at_points_RE8700}
\end{figure}

\begin{figure}[hbtp!]
\centering 
\includegraphics[width=\linewidth]{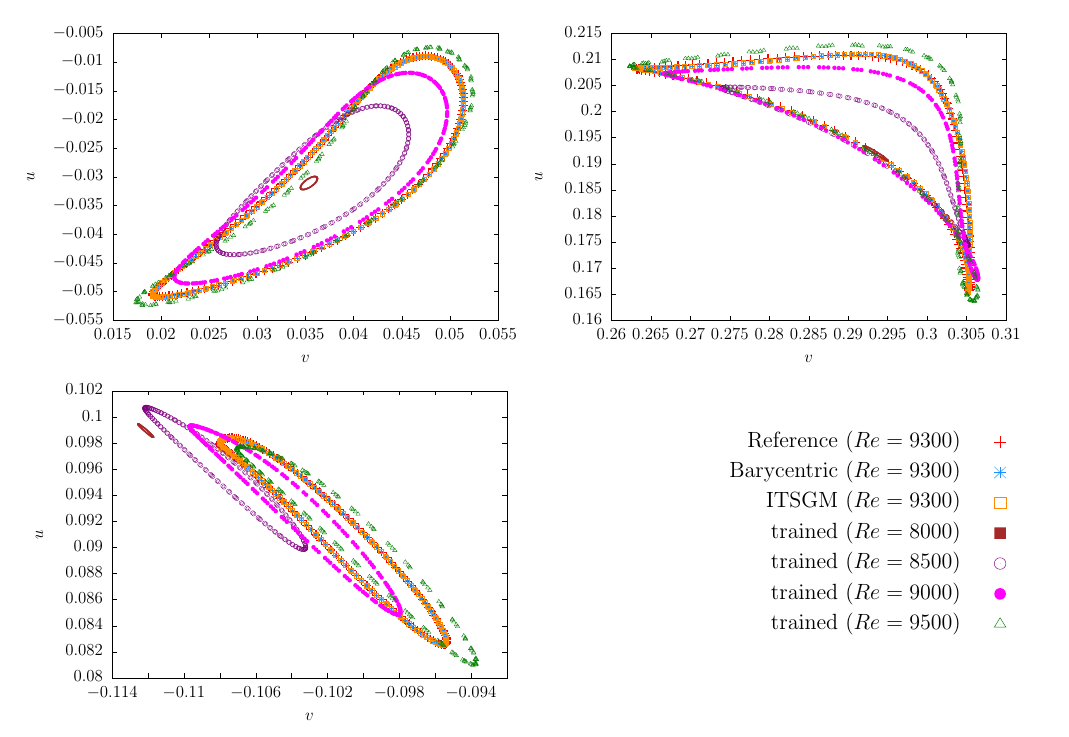}
\caption{Velocity phase portrait at the \Mourad{points $x_1$ (top left), $x_2$ (top right) and $x_3$ (bottom left)} for the untrained value $Re_1^{\txt{untr}} = 9300$. 
}
\label{fig:velocity_at_points_RE9300}
\end{figure}

Consider the two dimensional lid driven cavity flow where the space domain consists of a square cavity $\Omega = ]0,D[\times]0,D[$. At the top boundary, a tangential velocity $U$ of unit magnitude is applied to drive the fluid flow in the cavity, while no-slip conditions are set on the remaining boundaries. The characteristic length and velocity are set respectively to $D=1$ and $U = 1$.
The time step is set to $\delta t = 0.001$ and a non uniform mesh containing $32928$ elements is used for the triangulation of the cavity domain. A representation of the flow  dynamics at $Re=9000$ in three different time instants is given in figure \ref{fig:dynamics_cavity_flow}. A rich dynamics of the flow represented by the development of vortices in the corners can be noticed.
\vspace{0.2cm}
\\
A database of flow solutions is constructed by precomputing the solutions of the cavity flow problem in four trained Reynolds number values. These  values are $Re_1^{\txt{tr}}=8000$, $Re_2^{\txt{tr}}=8500$, $Re_3^{\txt{tr}}=9000$ and $Re_4^{\txt{tr}}=9500$. 
%
%
%A sketch up of the velocity fields highlighting the rich dynamics of the flow is depicted in figure \ref{fig:Comparison_sampling_sol_cavity}.
%
%
%
% 
%
%
%
%
%
For each value $Re_i^{\txt{tr}}$, $200$ snapshots covering nearly $10$ periods of the established regime of the flow are used to calculate the velocity POD bases of dimension $7$ each. In this study case, given that the POD vector functions are zero on the boundary and that there is no interest in calculating the pressure, the PROM \eqref{EQ.ROM_Navier_Stokes} describing only the velocity field is used. 
To this end, the offline stage of algorithm \ref{Alg:PROM_NS} is performed once for all. 
The aim of this application is to predict by interpolation, the solutions for the new untrained values $Re_1^{\txt{untr}}=8300$, $Re_2^{\txt{untr}}=8700$ and $Re_3^{\txt{untr}}=9300$. 
\vspace{0.2cm}
\\ 
Figures  \ref{fig:Mean_Errors} and \ref{fig:Errors_time} reports the mean error and its time history obtained by using both PROM approaches (Barycentric and ITSGM PROMs). In terms of accuracy, it can be seen that for the three test cases, both the proposed barycentric PROM and ITSGM PROM lead to nearly good equivalent errors (less than $0.6\%$). The results agreement between the two methods can be further inspected closely from the flow phase portraits at the three spatial points $x_1 = (2/16, 13/16)$, $x_2 = (2/16, 2/16)$ and $x_3 = (19/20, 19/20)$. These points are chosen in the corners of the cavity where the dynamics is rich. 
%The phase portraits are graphically reported in figures \ref{fig:velocity_at_points_RE8300}, \ref{fig:velocity_at_points_RE8700} and \ref{fig:velocity_at_points_RE9300}. 
Figures \ref{fig:velocity_at_points_RE8300}, \ref{fig:velocity_at_points_RE8700} and \ref{fig:velocity_at_points_RE9300} reveal a very good tracking of the high fidelity punctual flow trajectory either by Barycentric PROM  or by the ITSGM PROM.
%It can be observed that the Barycentric PROM reproduces solutions that track well and equivalently as the ITSGM PROM the punctual flow history produced by the High fidelity solution. 
%Moreover, from a global point of view, the reconstructed solutions at three instants of the flow depicted in figures \ref{fig:Comparison_sol_cavity_RE8300}, \ref{fig:Comparison_sol_cavity_RE8700} and \ref{fig:Comparison_sol_cavity_RE9300}, reveal the very good agreement of the new Barycentric PROM and the ITSGM PROM with the high fidelity solutions.
%
%\vspace{0.2cm}
%\\
%So far, It is numerically proven that the barycentric PROM and ITSGM PROM are nearly equivalent in terms of accuracy. 
The potential of the barycentric PROM over the ITSGM PROM appears  in the stage of construction of the reduced projection matrices. The time needed to update the barycentric PROM is of $1.1\times 10^{-3}$ seconds, while it is of $2.7s$ by the ITSGM PROM. In this case also, almost three orders of magnitude of computational time in the online time are gained by the new barycentric approach. This result supports the statement made about the effectiveness and flexibility of the proposed PROM approach in real-time applications. 
%This makes our approach very competitive for near real time applications.  

%\begin{figure}[hbtp!]
%\resizebox{0.5\textwidth}{!}{\large \input{./Figures/Errors/Mean_Error_velocity}}%
%\caption{Percentage of the mean errors of the solutions by the barycentric and ITSGM PROMs.}
%\label{fig:Mean_Errors}
%\end{figure}

%\begin{figure}[hbtp!]
%\hspace*{-2cm}
%\begin{subfigure}{0.5\textwidth}
%\resizebox{\textwidth}{!}{\large \input{./Figures/Errors/Error_velocity_RE8300}}%
%\caption{Errors for $Re_1^{\txt{untr}} = 8300$.}
%\end{subfigure}%
%\begin{subfigure}{0.5\textwidth}
%\resizebox{\textwidth}{!}{\large \input{./Figures/Errors/Error_velocity_RE8700}}
%\caption{Errors for $Re_2^{\txt{untr}} = 8700$.}
%\end{subfigure}
%\hspace*{-0.9cm}
%\begin{subfigure}{0.5\textwidth}
%\resizebox{\textwidth}{!}{\large \input{./Figures/Errors/Error_velocity_RE9300}}%
%\caption{Errors for $Re_3^{\txt{untr}} = 8300$.}
%\end{subfigure}%
%%\centering
%\hspace*{1.6cm}
%\hspace*{-1cm}\adjustbox{width= 0.55\linewidth,trim= 1cm 4cm 1cm 4cm}{\large \input{./Figures/Errors/LEGEND_Errors}}
%\caption{Time history of the percentage of error of the solutions obtained by the barycentric and ITSGM PROMs.}
%\label{fig:Errors_time}
%\end{figure}

%
%\begin{figure}[hbtp!]
%\hspace*{-2cm}
%\begin{subfigure}{0.5\textwidth}
%\resizebox{\textwidth}{!}{\large \input{./Figures/velocities_at_points/velocity_RE8300__p1}}%
%%\caption{Point $x_1 = (2/16, 13/16)$.}
%\end{subfigure}%
%\begin{subfigure}{0.5\textwidth}
%\resizebox{\textwidth}{!}{\large \input{./Figures/velocities_at_points/velocity_RE8300__p2}}%
%%\caption{Point $x_2 = (2/16, 2/16)$.}
%\end{subfigure}
%\hspace*{-0.9cm}
%\begin{subfigure}{0.5\textwidth}
%\resizebox{\textwidth}{!}{\large \input{./Figures/velocities_at_points/velocity_RE8300__p3}}%
%%\caption{Point $x_3 = (19/20, 19/20)$.}
%\end{subfigure}%
%%\centering
%\hspace*{1.6cm}
%\hspace*{-1cm}\adjustbox{width= 0.55\linewidth,trim= 1cm 4cm 1cm 4cm}{\large \input{./Figures/velocities_at_points/LEGEND_portrait_phase_RE8300}}
%\caption{Velocity phase portrait at the \Mourad{points $x_1$ (top left), $x_2$ (top right) and $x_3$ (bottom left)} for the untrained value $Re_1^{\txt{untr}} = 8300$. 
%}
%\label{fig:velocity_at_points_RE8300}
%\end{figure}

%\begin{figure}[hbtp!]
%\hspace*{-2cm}
%\begin{subfigure}{0.5\textwidth}
%\resizebox{\textwidth}{!}{\large \input{./Figures/velocities_at_points/velocity_RE8700__p1}}%
%%\caption{Point $x_1 = (2/16, 13/16)$.}
%\end{subfigure}%
%\begin{subfigure}{0.5\textwidth}
%\resizebox{\textwidth}{!}{\large \input{./Figures/velocities_at_points/velocity_RE8700__p2}}%
%%\caption{Point $x_2 = (2/16, 2/16)$.}
%\end{subfigure}
%\hspace*{-0.9cm}
%\begin{subfigure}{0.5\textwidth}
%\resizebox{\textwidth}{!}{\large \input{./Figures/velocities_at_points/velocity_RE8700__p3}}%
%%\caption{Point $x_3 = (19/20, 19/20)$.}
%\end{subfigure}%
%%\centering
%\hspace*{1.6cm}
%\hspace*{-1cm}\adjustbox{width= 0.55\linewidth,trim= 1cm 4cm 1cm 4cm}{\large \input{./Figures/velocities_at_points/LEGEND_portrait_phase_RE8700}}
%\caption{Velocity phase portrait at \Mourad{points $x_1$ (top left), $x_2$ (top right) and $x_3$ (bottom left)} for the untrained value $Re_2^{\txt{untr}} = 8700$. 
%%Here, $u$ and $v$ are respectively the horizontal and vertical components of the flow velocity $\bm{u}$.
%}
%\label{fig:velocity_at_points_RE8700}
%\end{figure}

%\begin{figure}[hbtp!]
%\hspace*{-2cm}
%\begin{subfigure}{0.5\textwidth}
%\resizebox{\textwidth}{!}{\large \input{./Figures/velocities_at_points/velocity_RE9300__p1}}%
%%\caption{Point $x_1 = (2/16, 13/16)$.}
%\end{subfigure}%
%\begin{subfigure}{0.5\textwidth}
%\resizebox{\textwidth}{!}{\large \input{./Figures/velocities_at_points/velocity_RE9300__p2}}
%%\caption{Point $x_2 = (2/16, 2/16)$.}
%\end{subfigure}
%\hspace*{-0.9cm}
%\begin{subfigure}{0.5\textwidth}
%\resizebox{\textwidth}{!}{\large \input{./Figures/velocities_at_points/velocity_RE9300__p3}}%
%%\caption{Point $x_3 = (19/20, 19/20)$.}
%\end{subfigure}%
%%\centering
%\hspace*{1.6cm}
%\hspace*{-1cm}\adjustbox{width= 0.55\linewidth,trim= 1cm 4cm 1cm 4cm}{\large \input{./Figures/velocities_at_points/LEGEND_portrait_phase_RE9300}}
%\caption{Velocity phase portrait at \Mourad{points $x_1$ (top left), $x_2$ (top right) and $x_3$ (bottom left)} for the untrained value $Re_3^{\txt{untr}} = 9300$. %Here, $u$ and $v$ are respectively the horizontal and vertical components of the flow velocity $\bm{u}$.
%}
%\label{fig:velocity_at_points_RE9300}
%\end{figure}
%
%
\section{Conclusions}
The present paper presents a new barycentric PROM approach based on the geometry of the quotient manifold $\mbb{R}^{N\times \nbrModes}_*/\mcal{O}(q)$ formed by the quotient of the set of full-rank $N$-by-$q$ matrices by the orthogonal group of dimension $q$. The subspace of projection is sought after as the center of mass of some trained subspaces (spanned by parametric POD bases) which is formulated as a minimization problem, and eventually solved as a simple fixed point problem. It was numerically demonstrated on the problems of the flow past a circular cylinder and the flow in a lid driven cavity, that the new barycentric and ITSGM PROMs provide nearly equivalent results. However, the potential of the barycentric approach over ITSGM reveals in the natural construction of the barycentric PROM. This last enables significant time speedups (several orders of magnitude) and thus real-time update of the ROM coefficients. With these significantly lower online computational resources, the proposed barycentric PROM is considered as a promising technique for parametric problems such as design and optimal control.

%%%%%%%%%%%%%%%%%%%%%%%%%%%%%%%%%%%%%%%%%%%%%%%%%%%%%%%%%%%%%%%%%%%%%%%%%%%%%%%%%%%%%%%%%%%%%%%%%%%%%%%%%%%%%%%%%%%%%%%%%%%%%%%%%%%%%%%%%%%%%%%%%%%
\section*{Acknowledgement}
%The authors would like to thank the Nouvelle-Aquitaine region and the European union for their financial support. 
This material is based upon work financially supported by CPER BATIMENT DURABLE - Axe 3 "Qualité des Environnement Intérieurs (QEI)" (P-2017-BAFE-102) and French Astrid ANR MODULO'PI (ANR-16-ASTR-0018 MODUL’O $\Pi$).
\newpage
\appendix
.
\section{Proper Orthogonal Decomposition}
Let $\{u_1,\dots,u_{\nbrSnap} \}$ be an ensemble of snapshots obtained by solving a physical problem. The aim of the POD is to represent this ensemble in the low dimensional subspace $span(\Xbasis{})$ that captures the quasi-totality of the dynamics of the underlying physical problem. The POD method consists of the following steps
\begin{itemize}
\item[\hspace{10pt}]
\begin{itemize}
\item[\textit{step 1}] build the correlation matrix $C$ as $C_{ij} =  \integOmega{ u_i \, u_j}$
\item[\textit{step 2}] solve the eigenvalue problem $C \Tbasis{} = \Tbasis{} \lambda $
\item[\textit{step 3}] calculate the POD modes $\SmodeU{k}{} = \somme{i}{1}{\nbrSnap}\lambda_k^{-\frac{1}{2}}  \Tbasis{ik} u_i $
%\item[\textit{step 4}] calculate the POD right singular vectors as $\Tbasis{} = $
\end{itemize}
\end{itemize}
In practice, the POD basis is truncated to an order $\nbrModes<\nbrSnap$, where only the first modes corresponding the the first $\nbrModes$ significant eigenvalues of $C$ are considered. This truncated basis is sufficient to represent most information contained in the original ensemble of snapshots.
\section{ITSGM method}\label{Appendix_ITSGM}
Let $span(\Sbase{1}), span(\Sbase{2}), \dots, span(\Sbase{\nbrParam})$ be a set of subspaces obtained from the POD of a parametrized physical problem at the trained points $\nu_1,\nu_2,\dots,\nu_{\nbrSnap}$. The ITSGM method is a standard method that enable to predict the subspace $span(\widetilde{\Sbase{}})$ associated to a new untrained parameter $\tilde{\nu}\neq\nu_k$. The steps of ITSGM are summarized in algorithm \ref{Alg:ITSGM} (see \cite{Amsallem} for more details).
\\
\begin{algorithm}[H]
\begin{itemize}
\item[\hspace{10pt}]
\begin{itemize}
\item[ \textit{step 1}] Choose the origin point of tangency $span(\Sbase{k_0})$ where $k_0 \in \{1,\dots,\nbrParam\}$.
\item[ \textit{step 2}] For $k \in \{1,\dots,\nbrParam\}$, map via the logarithmic mapping the point $span(\Sbase{k}) \in \GrassmanifoldSpace{}$ to the tangent space $\TangsubSpaceGrass{\Xbasis{k_0}}{}$ and obtain
$$ \InitVel{k} = U_k \arctan(\Sigma_k) V_k^T $$
where $(I-\Sbase{k_0} (\Sbase{k_0})^T)\Sbase{i}((\Sbase{k_0})^T\Sbase{i})^{- 1}((\Sbase{k_0})^T\Sbase{k_0})^{-1/2} = U_k \Sigma_k V_k^T$ (thin SVD).
\item[ \textit{step 3}] Interpolate the initial velocities $\InitVel{1}, \InitVel{2}, \dots, \InitVel{\nbrParam}$ for the untrained parameter $\tilde{\nu}$ using a standard interpolation and obtain $\widetilde{\xi}$.
\item[ \textit{step 4}] Finally by the exponential mapping, map the interpolated velocity $\xi$ back to the Grassmann manifold. The matrix spanning interpolated subspace is given by
$$
    \widetilde{\Sbase{}} = \Sbase{k_0} ((\Sbase{k_0})^T\Sbase{k_0})^{-1/2} {V} \cos({\Sigma}) + {U} \sin({\Sigma})
$$
where $U \Sigma V^T$ is the thin SVD of the interpolated initial velocity vector $\widetilde{\xi}$.
\end{itemize}
\end{itemize}
\caption{ITSGM}
\label{Alg:ITSGM}
\end{algorithm}
\section{Velocity/pressure reduced order model}\label{Appendix_ROM2}
Let $\Smode{u}{}$ and $\Smode{p}{}$ be the velocity and pressure POD bases respectively of dimensions $\nbrModesU$ and $\nbrModesP$ of the problem of the flow past a cylinder. The reduced order model describing the the evolution of both velocity and pressure is obtained by the orthogonal projection of the residual onto the velocity basis functions and onto the gradient of the pressure basis functions \cite{Tallet-Allery-com-2015}. This ROM writes as follows
\begin{equation}\label{Eq : ROM CYLINDER}
\begin{cases}
\Mat{M}^{(u)} \derivee{\Tmode{u}{}}{t} +  \nu \Mat{R}^{(u)} \Tmode{u}{}+  \Mat{\overline{C}}^{(u)} \Tmode{u}{} +  \somme{e}{1}{\nbrModesU} \Tmode{u}{e}\Mat{C}^{(u),k} \Tmode{u}{} + \Mat{K}^{(u)} \Tmode{p}{}= {\bm{F}}^{(u)}, \ 
\\
\Mat{M}^{(p)} \derivee{\Tmode{u}{}}{t} +  \nu \Mat{R}^{(p)} \Tmode{u}{}+\Mat{\overline{C}}^{(p)}  \Tmode{u}{} + \somme{e}{1}{\nbrModesU} \Tmode{u}{e}\Mat{C}^{(p),j} \Tmode{u}{}+  \Mat{K}^{(p)} \Tmode{p}{}= {\bm{F}}^{(p)}
%\\
%\somme{j}{1}{\nbrModesU}\Mat{M}^{(u)}_{ij} \Tmode{u}{j}(0) = \int_{\Omega} \vect{u_0}  \Smode{u}{i}\,dx
%\\ 
%\forall i = 1,\cdots, \nbrModesU, \ \ \ \forall m = 1,\cdots, \nbrModesP
\end{cases}
\end{equation}
where 
\begin{align*}
\Mat{M}^{(u)}_{ij} &= \integOmega{\Smode{u}{j} \Smode{u}{i}}
\hspace*{1cm}
\Mat{R}^{(u)}_{ij}  = \integOmega{\nabla \Smode{u}{j} : \nabla\Smode{u}{i}} 
\hspace*{1cm} 
\Mat{K}^{(u)}_{il}  = \int_{\Gamma_{\txt{horiz}}} \Smode{p}{l} \Smode{u}{i}\cdot \vect{n} \,d\sigma
\\
\Mat{\overline{C}}^{(u)}_{ij} &= \integOmega{(\overline{\statevecU}\cdot\nabla) \Smode{u}{j} \cdot \Smode{u}{i}} + \integOmega{(\Smode{u}{j}\cdot\nabla) \overline{\statevecU} \cdot \Smode{u}{i}}\\
\Mat{C}^{(u),e}_{ij} &= \integOmega{(\Smode{u}{e}\cdot\nabla) \Smode{u}{j} \cdot \Smode{u}{i}} 
\hspace*{2cm}
{\bm{F}}^{(u)}_{i} = \integOmega{ \left(\nu \Delta \overline{\statevecU} - \overline{\statevecU} \cdot\nabla \overline{\statevecU} - \nabla \overline{p}\right) \, \Smode{u}{i}}
\\
\Mat{M}^{(p)}_{mj} &= \integOmega{\Smode{u}{j} \nabla\Smode{p}{m}}
\hspace*{1cm}
\Mat{R}^{(p)}_{mj}  = -\integOmega{\Delta \Smode{u}{j} : \nabla\Smode{p}{m}} 
\hspace*{1cm} 
\Mat{K}^{(p)}_{ml}  = \int_{\Omega} \nabla\Smode{p}{l} \nabla\Smode{p}{m} \,dx
\\
\Mat{\overline{C}}^{(u)}_{mj} &= \integOmega{(\overline{\statevecU}\cdot\nabla) \Smode{u}{j} \cdot \nabla\Smode{p}{m}} + \integOmega{(\Smode{u}{j}\cdot\nabla) \overline{\statevecU} \cdot \nabla\Smode{p}{m}}\\
\Mat{C}^{(p),e}_{mj} &= \integOmega{(\Smode{u}{e}\cdot\nabla) \Smode{u}{j} \cdot \nabla\Smode{p}{m}} 
\hspace*{2cm}
{\bm{F}}^{(p)}_{m} = \integOmega{ \left(\nu \Delta \overline{\statevecU} - \overline{\statevecU} \cdot\nabla \overline{\statevecU} - \nabla \overline{p}\right) \, \nabla\Smode{p}{m}}
\end{align*}

\section{Aerodynamics coefficients}\label{Appendix_aerodynamics_coeffs}
The Strouhal number is defined by 
$$S_t = \myfrac{D f_s}{U} $$ 
where $D$ is the characteristic length, $U$ the characteristic velocity and $f_s$ the oscillation frequency of the flow. The drag and lift coefficients $C_D$ and $C_L$ are calculated as follows
\begin{equation*}
C_D = \myfrac{F_D}{\frac{1}{2}\rho U^2 D} \ \ ; \ \  C_L = \myfrac{F_L}{\frac{1}{2}\rho U^2 D}
\end{equation*}
where $\rho$ is the density of the fluid and $F_D$, $F_L$ represent respectively the drag and lift forces exerted by the fluid on the cylinder. These forces are calculated as follows
\begin{equation*}
\left(\begin{matrix}
F_D \\
F_L
\end{matrix}\right)  =  \int_{\Gamma_{circ}} \left( \mu (\nabla \vect{u} + \nabla^T \vect{u}) - p I \right)\cdot \vect{n} \,d\sigma
\end{equation*} 
%
%\section{Some precisions on the Bi-CITSGM}\label{Appendix Bi-CITSGM}
 %%%%%%%%%%%%%%%%%%%%%%%%%%%%%%%%%%%%%%%%%%%%%%%%%%%%%%%%%%%%%%%%%%%%%%%%%%%%%%%%%%%%%%%%%%%%%%%%%%%%%%%%%%%%%%%%%%%%%%%%%%%%%%%%%%%%%%%%%%%%%%%%%%%%

\section*{\large References}
\bibliographystyle{ieeetr}
\bibliography{./BIBLIO}
\end{document}